\documentclass[a4, 10pt]{amsart}
\usepackage{amssymb}
\usepackage{amstext}
\usepackage{amsmath}
\usepackage{amscd}
\usepackage{latexsym}
\usepackage{amsfonts}

\theoremstyle{plain}
\newtheorem{thm}{Theorem}[subsection]
\newtheorem*{thma}{Theorem A}
\newtheorem*{thmb}{Theorem B}

\newtheorem{prop}[thm]{Proposition}
\newtheorem{lem}[thm]{Lemma}

\newtheorem{cor}[thm]{Corollary}
\newtheorem{claim}{Claim}

\theoremstyle{definition}
\newtheorem{defn}[thm]{Definition}

\newtheorem{ques}[thm]{Question}

\newtheorem*{ac}{\sc Acknowledgments}

\theoremstyle{remark}
\newtheorem*{pf}{{\sc Proof}}

\newtheorem*{tpf}{{\sc Proof of Theorem \ref{2dimgor}}}

\numberwithin{equation}{thm}
\def\Hom{\mathrm{Hom}}

\def\Ext{\mathrm{Ext}}

\def\End{\mathrm{End}}
\def\Soc{\mathrm{Soc}}

\def\mod{\mathrm{mod\,}}
\def\Coker{\mathrm{Coker\,}}
\def\Ker{\mathrm{Ker\,}}
\def\Im{\mathrm{Im\,}}

\def\gr{\mathrm{gr}}

\def\a{\mathfrak a}
\def\b{\mathfrak b}
\def\m{\mathfrak m}
\def\n{\mathfrak n}

\def\g{\mathrm{G}\textrm{-}\mathrm{dim}}
\def\Kdim{\mathrm{dim}}
\def\edim{\mathrm{edim\,}}
\def\depth{\mathrm{depth}}

\def\Ann{\mathrm{Ann}}

\def\pd{\mathrm{pd}}

\def\height{\mathrm{ht\,}}
\def\grade{\mathrm{grade}}

\def\C{{\mathcal C}}
\def\F{{\mathcal F}}

\def\G{{\mathcal G}}

\begin{document}

%%%%%%%%%%%%%%%%%%% Title %%%%%%%%%%%%%%%%%%%%%%%%%%%%%%%%%
\title{Direct summands of syzygy modules of the residue class field}
\author{Ryo Takahashi}
\address{Department of Mathematics, School of Science and Technology, Meiji University, Kawasaki 214-8571, Japan}
\email{takahasi@math.meiji.ac.jp}
%%%%%%%%%%%%%%%%%%% Abstract %%%%%%%%%%%%%%%%%%%%%%%%%%%%%
\begin{abstract}
Let $R$ be a commutative Noetherian local ring.
This paper deals with the problem asking whether $R$ is Gorenstein if the $n$th syzygy module of the residue class field of $R$ has a non-trivial direct summand of finite G-dimension for some $n$.
It is proved that if $n$ is at most two then it is true, and moreover, the structure of the ring $R$ is determined essentially uniquely.
\end{abstract}
\maketitle
%%%%%%%%%%%%%%%%%%%%%%%%%%%%%%%%%%%%%%%%%%%%%%%%%%%%%%%%%
%\tableofcontents
%%%%%%%%%%%%%%%%%%% Section 1 %%%%%%%%%%%%%%%%%%%%%%%%%%
\section{Introduction}

Throughout the present paper, we assume that all rings are commutative Noetherian local rings and all modules are finitely generated modules.

G-dimension is a homological invariant of a module which has been introduced by Auslander \cite{Auslander}.
This invariant is an analogue of projective dimension.
Whereas the finiteness of projective dimension characterizes the regular property of the base ring, the finiteness of G-dimension characterizes the Gorenstein property of the base ring.
To be precise, any module over a Gorenstein local ring has finite G-dimension, and a local ring with residue class field of finite G-dimension is Gorenstein.
G-dimension shares a lot of properties with projective dimension.
For example, it also satisfies an Auslander-Buchsbaum-type equality, which is called the Auslander-Bridger formula.

Dutta \cite{Dutta} proved the following theorem in his research into the homological conjectures:

\begin{thm}[Dutta]\label{dutta}
Let $(R, \m, k)$ be a local ring.
Suppose that the $n$th syzygy module of $k$ has a non-zero direct summand of finite projective dimension for some $n \geq 0$.
Then $R$ is regular.
\end{thm}

Since G-dimension is similar to projective dimension, this theorem naturally leads us to the following question:

\begin{ques}\label{ques}
Let $(R, \m, k)$ be a local ring.
Suppose that the $n$th syzygy module of $k$ has a non-zero direct summand of finite G-dimension for some $n \geq 0$.
Then is $R$ Gorenstein?
\end{ques}

It is obviously seen from the indecomposability of $k$ that this question is true if $n=0$.
Hence this question is worth considering just in the case where $n\geq 1$.

We are able to answer in this paper that the above question is true if $n\leq 2$.
Furthermore, as the theorems below say, we can even determine the structure of a ring satisfying the assumption of the above question for $n=1, 2$.

The organization of this paper is as follows.
In Section 2, we will prepare some notions and results for later use.
The definition and properties of G-dimension will be given in this section.
In Section 3, we shall state the main theorems of this paper.
Firstly, we will consider a local ring such that the first syzygy module of the residue class field, namely, the maximal ideal, is decomposable.
We will obtain the following result:

\begin{thma}
Let $(R, \m )$ be a complete local ring.
The following conditions are equivalent:
\begin{enumerate}
\item[$(1)$]\ There is an $R$-module $M$ with $\g _R M < \infty = \pd _R M$, and $\m$ is decomposable;
\item[$(2)$]\ $R$ is Gorenstein, and $\m$ is decomposable;
\item[$(3)$]\ There are a complete regular local ring $S$ of dimension two and a regular system of parameters $x, y$ of $S$ such that $R\cong S/(xy)$.
\end{enumerate}
\end{thma}

Secondly, we will investigate a local ring such that the second syzygy module of the residue class field is decomposable, and obtain the following result:

\begin{thmb}
Let $(R, \m , k)$ be a complete local ring.
{\rm (Denote by $\Omega _R ^2 k$ the second syzygy module of $k$.)}
Suppose that $\m$ is indecomposable.
Then the following conditions are equivalent:
\begin{enumerate}
\item[$(1)$]\ There is a non-trivial direct summand $M$ of $\Omega _R ^2 k$ with $\g _R M<\infty$;
\item[$(2)$]\ $R$ is Gorenstein, and $\Omega _R ^2 k$ is decomposable;
\item[$(3)$]\ There are a complete regular local ring $(S, \n )$ of dimension three, a regular system of parameters $x, y, z$ of $S$, and $f\in\n$ such that $R\cong S/(xy-zf)$.
\end{enumerate}
\end{thmb}

Theorem A and B especially say that a complete Gorenstein local ring such that the first or second syzygy module of the residue class field is decomposable is a hypersurface, and moreover, its ring structure can be determined concretely.
We will actually prove in Section 3 more general results than the above two theorems.

%%%%%%%%%%%%%%%%%%% Section 2 %%%%%%%%%%%%%%%%%%%%%%%%%%%%%
\section{Preliminaries}

Throughout this section, let $(R, \m , k)$ be a local ring.
In this section, we will recall several basic notions and state related results to explain and prove the main theorems of this paper.

%%%%% (Pre)covers and (pre)envelopes %%%%%
\subsection{(Pre)covers and (pre)envelopes}

We begin by recalling the notions of a (pre)cover and a (pre)envelope of a module.
Let $\mod R$ denote the category of finitely generated $R$-modules.

\begin{defn}
Let $\C$ be a full subcategory of $\mod R$.\\
(1) Let $\phi : X\to M$ be a homomorphism from $X\in\C$ to $M\in\mod R$.

(i) We call $\phi$ or $X$ a {\it $\C$-precover} of $M$ if for any homomorphism $\phi ': X'\to M$ with $X'\in\C$ there exists a homomorphism $f: X'\to X$ such that $\phi '=\phi f$.

(ii) Assume that $\phi$ is a $\C$-precover of $M$.
We call $\phi$ or $X$ a {\it $\C$-cover} of $M$ if any endomorphism $f$ of $X$ with $\phi =\phi f$ is an automorphism.\\
(2) Let $\phi : M\to X$ be a homomorphism from $M\in\mod R$ to $X\in\C$.

(i) We call $\phi$ or $X$ a {\it $\C$-preenvelope} of $M$ if for any homomorphism $\phi ': M\to X'$ with $X'\in\C$ there exists a homomorphism $f: X\to X'$ such that $\phi '=f\phi$.

(ii) Assume that $\phi$ is a $\C$-preenvelope of $M$.
We call $\phi$ or $X$ a {\it $\C$-envelope} of $M$ if any endomorphism $f$ of $X$ with $\phi =f\phi$ is an automorphism.
\end{defn}

A $\C$-precover (resp. $\C$-cover, $\C$-preenvelope, $\C$-envelope) is also called a {\it right $\C$-approximation} (resp. {\it right minimal $\C$-approximation}, {\it left $\C$-approximation}, {\it left minimal $\C$-approximation}).
A $\C$-cover (resp. $\C$-envelope) is uniquely determined up to isomorphism whenever it exists.
In general, it is uncertain whether the existence of a $\C$-precover (resp. $\C$-preenvelope) implies the existence of a $\C$-cover (resp. $\C$-envelope). 
However, it is true under a few assumptions: if $R$ is Henselian and $\C$ is closed under direct summands, then for a given $\C$-precover (resp. $\C$-preenvelope), we can extract a $\C$-cover (resp. $\C$-envelope) from it, as follows.

\begin{prop}
Let $\C$ be a full subcategory of $\mod R$ which is closed under direct summands.
Suppose that $R$ is Henselian.
\begin{enumerate}
\item[$(1)$]\ Let $0 \to N \to X \overset{\phi}{\to} M$ be an exact sequence of $R$-modules where $\phi$ is a $\C$-precover of $M$.
Then there exists a commutative diagram
$$
\begin{CD}
@. 0 @. 0 \\
@. @VVV @VVV \\
@. L @= L \\
@. @VVV @VVV \\
0 @>>> N @>>> X @>{\phi}>> M \\
@. @VVV @VVV @| \\
0 @>>> N' @>>> X' @>{\phi '}>> M \\
@. @VVV @VVV \\
@. 0 @. 0
\end{CD}
$$
of $R$-modules with exact rows and split exact columns such that $\phi '$ is a $\C$-cover of $M$.
\item[$(2)$]\ Let $M \overset{\phi}{\to} X \to N \to 0$ be an exact sequence of $R$-modules where $\phi$ is a $\C$-preenvelope of $N$.
Then there exists a commutative diagram
$$
\begin{CD}
@. 0 @. 0 \\
@. @VVV @VVV \\
M @>{\phi '}>> X' @>>> N' @>>> 0 \\
@| @VVV @VVV \\
M @>{\phi}>> X @>>> N @>>> 0 \\
@. @VVV @VVV \\
@. L @= L \\
@. @VVV @VVV \\
@. 0 @. 0
\end{CD}
$$
of $R$-modules with exact rows and split exact columns such that $\phi '$ is a $\C$-envelope of $M$.
\end{enumerate}
\end{prop}

For the proof of the statement (1), we refer to \cite[Remark 2.6]{Takahashi}.
The statement (2) can be shown dually.

%%%%% The subcategory of free modules %%%%%
\subsection{The subcategory of free modules}

We denote by $\F (R)$ the full subcategory of $\mod R$ consisting of all free $R$-modules.
Recall that a homomorphism $f : M \to N$ of $R$-modules is said to be {\it minimal} if the induced homomorphism $f\otimes _R k : M\otimes _R k \to N\otimes _R k$ is an isomorphism.
(Note from Nakayama's lemma that every minimal homomorphism is surjective.)
Let $\nu _R (M)$ denote the minimal number of generators of an $R$-module $M$, i.e., $\nu _R (M)=\Kdim _k (M\otimes _R k)$.
Set $(-)^{\ast}=\Hom _R (-, R)$.
Every $R$-module admits an $\F(R)$-cover and an $\F(R)$-envelope, as follows.

\begin{prop}
Let $M$ be an $R$-module.
\begin{enumerate}
\item[{\rm (1)}]
A homomorphism $\phi : R^n \to M$ is an $\F (R)$-cover of $M$ if and only if $\phi$ is surjective and $n=\nu _R (M)$.
\item[{\rm (2)}]
Let $f_1, f_2, \dots, f_n$ be a minimal system of generators of $M^{\ast}$.
Then the homomorphism $f={}^t(f_1,\dots,f_n): M \to R^n$ is an $\F (R)$-envelope of $M$.
\end{enumerate}
\end{prop}

An $R$-module $M$ is said to be {\it torsionless} (resp. {\it reflexive}) if the natural homomorphism $M \to M^{\ast\ast}$ is injective (resp. bijective).
We easily obtain the following.

\begin{cor}
Let $M$ be an $R$-module.
\begin{enumerate}
\item[$(1)$]\ Let $\sigma : M \to M^{\ast\ast}$ be the natural homomorphism and $\phi : F \to M^{\ast}$ an $\F (R)$-cover.
Then the composite map $\phi ^{\ast} \sigma: M \to F^{\ast}$ is an $\F (R)$-envelope.
\item[$(2)$]\ The $R$-module $M$ is torsionless if and only if the $\F (R)$-envelope of $M$ is an injective homomorphism.
\end{enumerate}
\end{cor}

We especially see from this corollary that an $\F (R)$-envelope is not necessarily an injective homomorphism.

Let $M$ be an $R$-module.
Take its $\F (R)$-cover $\pi : F \to M$.
The {\it first syzygy module} $\Omega _R M=\Omega _R ^1 M$ of $M$ is defined to be the kernel of the homomorphism $\pi$, and the {\it $n$th syzygy module} $\Omega _R ^n M$ of $M$ is defined inductively: $\Omega _R ^n M=\Omega _R (\Omega _R ^{n-1} M)$ for $n\geq 2$.
Dually to this, we can define the cosyzygy modules of any module.

\begin{defn}
Let $M$ be an $R$-module.\\
(1) Take the $\F (R)$-envelope $\theta : M \to F$ of $M$.
We set $\Omega _R ^{-1} M = \Coker\theta$, and call it the {\it first cosyzygy module} of $M$.\\
(2) Let $n\geq 2$.
Assume that the $(n-1)$th cosyzygy module $\Omega _R ^{-(n-1)} M$ is defined.
Then we set $\Omega _R ^{-n} M = \Omega _R ^{-1} (\Omega _R ^{-(n-1)} M)$ and call it the {\it $n$th cosyzygy module} of $M$.
\end{defn}

A module is said to be {\it stable} if it has no non-zero free summand.
The following is a property which is peculiar to cosyzygy modules.

\begin{prop}
For any $M\in\mod R$ and any $n\geq 1$, the module $\Omega _R ^{-n} M$ is stable.
\end{prop}

%%%%% G-dimension %%%%%
\subsection{G-dimension}

Now, we recall the definition of G-dimension.

\begin{defn}
(1) We denote by $\G (R)$ the full subcategory of $\mod R$ consisting of all $R$-modules $M$ satisfying the following three conditions:

(i) $M$ is reflexive,

(ii) $\Ext _R ^i (M, R)=0$ for every $i>0$,

(iii) $\Ext _R ^i (M^{\ast}, R)=0$ for every $i>0$.\\
(2) Let $M$ be an $R$-module.
If $n$ is a non-negative integer such that there is an exact sequence $0 \to G_n \to G_{n-1} \to \cdots \to G_1 \to G_0 \to M \to 0$ of $R$-modules with $G_i\in\G (R)$ for every $i$, then we say that $M$ has {\it G-dimension at most $n$}, and write $\g _R M\leq n$.
If such an integer $n$ does not exist, then we say that $M$ has {\it infinite G-dimension}, and write $\g _R M=\infty$.
\end{defn}

If an $R$-module $M$ has G-dimension at most $n$ but does not have G-dimension at most $n-1$, then we say that $M$ has {\it G-dimension $n$}, and write $\g _R M=n$.
Note that for an $R$-module $M$ we have $\g _R M = 0$ if and only if $M\in\G (R)$, and that all free $R$-modules belong to $\G (R)$.
For basic properties of G-dimension, we should refer to \cite[Chapter 3,4]{AB}, \cite[Chapter 1]{Christensen} and \cite[Section 8]{Avramov2}.
We write down some properties of the category $\G (R)$.

\begin{prop}
\begin{enumerate}
\item[$(1)$]\ If $R$ is a Gorenstein local ring, then the category $\G (R)$ coincides with the full subcategory of $\mod R$ consisting of all maximal Cohen-Macaulay modules.
\item[$(2)$]\ There exists a non-free $R$-module in $\G (R)$ if and only if there exists an $R$-module of finite G-dimension and infinite projective dimension.
\item[$(3)$]\ The following statements hold:
\begin{itemize}
\item[{\rm (i)}]\ If an $R$-module $M$ belongs to $\G (R)$, then so do $M^{\ast}$, $\Omega M$, $\Omega ^{-1} M$;
\item[{\rm (ii)}]\ Let $M, N$ be $R$-modules.
Then $M, N$ belong to $\G (R)$ if and only if so does $M\oplus N$;
\item[{\rm (iii)}]\ Let $0 \to L \to M \to N \to 0$ be an exact sequence of $R$-modules.
If $L, N$ belong to $\G (R)$, then so does $M$;
\item[{\rm (iv)}]\  If an $R$-module $M$ belongs to $\G (R)$, the $R/(x)$-module $M/xM$ belongs to $\G (R/(x))$ for any element $x\in\m$ which is $R$- and $M$-regular.
\end{itemize}
\end{enumerate}
\end{prop}

If $R$ is Gorenstein and non-regular, then the latter condition in (2) of the above proposition holds.
In fact, the $R$-module $k$ has finite G-dimension and infinite projective dimension.

We denote by $\underline{\G}(R)$ the full subcategory of $\G (R)$ consisting of all stable modules in $\G (R)$.
The dual functor $(-)^{\ast}$ and the syzygy functor $\Omega (-)$ make good correspondences between the category $\underline{\G}(R)$ and itself.

\begin{prop}
\begin{enumerate}
\item[$(1)$]\ We have an anti-equivalence of categories
$$
\underline{\G}(R)\to\underline{\G}(R),\ M\mapsto M^{\ast}
$$
with the functor being its own quasi-inverse.
\item[$(2)$]\ We have an equivalence of categories
$$
\underline{\G}(R)\to\underline{\G}(R),\ M\mapsto\Omega M
$$
having as quasi-inverse functor the functor $\underline{\G}(R)\to\underline{\G}(R),\ M\mapsto\Omega ^{-1}M$.
\end{enumerate}
\end{prop}

This proposition yields the following corollary.

\begin{cor}
For an $R$-module $M$, the following are equivalent:
\begin{enumerate}
\item[$(1)$]\ $M$ is a non-free indecomposable module in $\G (R)$;
\item[$(2)$]\ $M^{\ast}$ is a non-free indecomposable module in $\G (R)$;
\item[$(3)$]\ $\Omega M$ is a non-free indecomposable module in $\G (R)$;
\item[$(4)$]\ $\Omega ^{-1}M$ is a non-free indecomposable module in $\G (R)$.
\end{enumerate}
\end{cor}

%%%%% The fundamental module %%%%%
\subsection{The fundamental module}
Here we introduce the concept of the fundamental module.

\begin{defn}
Let $(R, \m , k)$ be a Cohen-Macaulay local ring of dimension two with canonical module $K$.
Then since $\Ext _R ^1 (\m , K) \cong \Ext _R ^2 (k, K) \cong k$, there exists a non-split exact sequence $\sigma : 0 \to K \to E \to \m \to 0$ which is unique up to equivalence.
This sequence $\sigma$ is called the {\it fundamental sequence} of $R$ and the intermediate module $E$ is called the {\it fundamental module} of $R$.
\end{defn}

We recall a numerical invariant of a module, which was invented by Auslander.

\begin{defn}
Let $R$ be a Gorenstein local ring.
For an $R$-module $M$, we denote by $\delta _R (M)$ the maximal rank of free summands of the $\G (R)$-cover of $M$, and set $\delta _R ^i (M)=\delta _R (\Omega _R ^i M)$, which is called the {\it $i$th Auslander's $\delta$-invariant} of $M$.
\end{defn}

\begin{lem}[Auslander]\label{auslander}
Let $R$ be a Gorenstein non-regular local ring with residue class field $k$.
Then $\delta _R ^i (k) = 0$ for every $i\geq 0$.
In other words, every syzygy module of $k$ admits a stable $\G (R)$-cover.
\end{lem}

This lemma was proved by Auslander in the unpublished paper \cite{Auslander2}.
For the proof, we can refer to \cite[Theorem 6]{Martsinkovsky}, \cite[Proposition 5.7]{ADS}, or \cite[Theorem (4.8)]{Yoshino2}.

Now, we can investigate several properties of the fundamental module of a Gorenstein local ring of dimension two.

\begin{prop}\label{fund}
Let $R$ be a Henselian Gorenstein non-regular local ring of dimension two, and let $\sigma : 0 \to R \to E \overset{\phi}{\to} \m \to 0$ be the fundamental sequence of $R$.
Then
\begin{enumerate}
\item[$(1)$]\ $\phi$ is the $\G (R)$-cover of $\m$,
\item[$(2)$]\ $E$ is stable,
\item[$(3)$]\ $E \cong \Omega _R ^{-1} (\Omega _R ^2 k)$,
\item[$(4)$]\ $E$ is indecomposable if and only if so is $\Omega _R ^2 k$.
\end{enumerate}
\end{prop}

\begin{pf}
(1) Since $R$ is Gorenstein, $\G(R)$ coincides with the category of maximal Cohen-Macaulay $R$-modules, and the assertion is a well-known fact on the fundamental sequence.

(2) This assertion follows from (1) and Lemma \ref{auslander}.

(3) Set $M=\Omega _R ^2 k$.
Note that the module $M$ belongs to $\G (R)$.
There is also an exact sequence $0 \to M \overset{\alpha}{\to} R^e \to \m \to 0$.
Take a minimal homomorphism $\beta : R^r \to M^{\ast}$.
Then $\Coker (\beta ^{\ast})$ is isomorphic to $\Omega _R ^{-1} M$ by definition.
The dual homomorphism $\alpha ^{\ast}: (R^e)^{\ast} \to M^{\ast}$ factors through $\beta$, i.e., there exists a homomorphism $\gamma : (R^e)^{\ast} \to R^r$ such that $\alpha ^{\ast}=\beta \gamma$.
Hence we have $\alpha = \alpha ^{\ast\ast} = \gamma ^{\ast} \beta ^{\ast}$.
Since $R\cong R^{\ast}$ and $M\cong M^{\ast\ast}$, we see that there is a commutative diagram
$$
\begin{CD}
0 @>>> M @>{\beta ^{\ast}}>> R^r @>{\varepsilon}>> \Omega _R ^{-1} M @>>> 0\\
@. @| @V{\gamma ^{\ast}}VV @VVV @.\\
0 @>>> M @>{\alpha}>> R^e @>>> \m @>>> 0
\end{CD}
$$
with exact rows.
Take a minimal homomorphism $\zeta : R^s \to \Coker (\gamma ^{\ast})$ and let $\eta : R^e \to \Coker (\gamma ^{\ast})$ be the natural surjection.
Then there is a homomorphism $\theta : R^s \to R^e$ such that $\zeta = \eta\theta$.
We easily see that the homomorphism $(\gamma ^{\ast} , \theta) : R^r \oplus R^s \to R^e$ is surjective, and obtain a commutative diagram
$$
\begin{CD}
@. @. 0 @. 0\\
@. @. @VVV @VVV\\
@. @. R^t @= R^t\\
@. @. @VVV @VVV\\
0 @>>> M @>{\left(
\begin{smallmatrix}
\beta ^{\ast}\\
0
\end{smallmatrix}
\right)}>> R^r \oplus R^s @>{\left(
\begin{smallmatrix}
\varepsilon & 0\\
0 & 1
\end{smallmatrix}
\right)}>> \Omega _R ^{-1} M \oplus R^s @>>> 0\\
@. @| @V{\left(
\begin{smallmatrix}
\gamma ^{\ast} & \theta
\end{smallmatrix}
\right)}VV @V{\kappa}VV\\
0 @>>> M @>{\alpha}>> R^e @>>> \m @>>> 0\\
@. @. @VVV @VVV\\
@. @. 0 @. 0
\end{CD}
$$
with exact rows and columns, where $t=r+s-e$.
Hence the homomorphism $\kappa$ is a $\G (R)$-precover of $\m$.
It follows from (1) that there is an isomorphism $\Omega _R ^{-1} M \oplus R^s \cong E\oplus R^{t-1}$.
Since both $E$ and $\Omega _R ^{-1} M$ are stable by (2), we conclude from the Krull-Schmidt theorem that the module $E$ is isomorphic to $\Omega _R ^{-1} M$, as desired.

(4) This assertion is proved by (3).
\qed
\end{pf}

%%%%%%%%%%% Section 3 %%%%%%%%%%%%%%%%%%
\section{Main results}

In this section, using the results given in the previous section, we shall state and prove our main theorems.

%%%%% Idealizations %%%%%
\subsection{Idealizations}

First of all, we consider an idealization possessing a non-free reflexive module.
We begin with making an easy lemma, which will often be used later.

\begin{lem}\label{nfi}
Let $(R, \m)$ be a local ring, $\theta : \m \to R$ the natural inclusion map, and $M$ a stable $R$-module.
Then the induced injective homomorphism
$$
\Hom _R (M, \theta ) : \Hom _R (M, \m ) \to M^{\ast}
$$
is an isomorphism.
\end{lem}

\begin{pf}
If there is a homomorphism from $M$ to $R$ which does not factor through $\theta$, then it is a surjection, hence is a split-epimorphism, contrary to the stability of $M$.
\qed
\end{pf}

Now we can prove the following result.

\begin{prop}\label{ideal}
Let $(S, \n , k)$ be a local ring, $V \not = 0$ a finite-dimensional $k$-vector space, and $R = S \ltimes V$ the idealization of $V$ over $S$.
Let $M$ be a non-free indecomposable reflexive $R$-module.
Then
\begin{enumerate}
\item[$(1)$]\ $M\cong\Soc\, R \cong V \cong k$,
\item[$(2)$]\ If $\depth\, S=0$, then $S=k$, hence $R\cong k[[X]]/(X^2)$.
\end{enumerate}
\end{prop}

\begin{pf}
(1)\ Denote by $\m$ the unique maximal ideal of $R$, and set $I =\n \ltimes 0 =\{(s,v)\in R \mid s\in\n, v=0\}$, and $J =0 \ltimes V =\{(s,v)\in R \mid s=0\}$.
These are ideals of $R$, and it is easy to see that $\m = I \oplus J$.
By virtue of Lemma \ref{nfi}, we have isomorphisms $M^{\ast} \cong \Hom _R (M, \m)\cong\Hom _R (M, I \oplus J)\cong\Hom _R (M, I)\oplus\Hom _R (M, J)$.
Since $M^{\ast}$ is also indecomposable, we have either $\Hom _R (M, I)=0$ or $\Hom _R (M, J)=0$.
However $J$ is isomorphic to $k^e$ as an $R$-module where $e=\Kdim _k V$, hence $\Hom _R (M, J)\cong k^{ne}\not =0$ where $n=\nu _R (M)$.
It follows that
\begin{equation}\label{mi}
\Hom _R (M, I)=0
\end{equation}
and $M^{\ast}\cong k^{ne}$.
The indecomposability of $M^{\ast}$ again implies that $M^{\ast}\cong k$ and $ne=1$, hence $e=1$.
Therefore $V\cong k$.
Also, we have isomorphisms $M\cong M^{\ast\ast}\cong k^{\ast}\cong k^r$ where $r=\Kdim _k (\Soc\, R)$.
The indecomposability of $M$ implies that $M\cong k$ and $r=1$.
Hence $\Soc\, R\cong k$.

(2)\ Note from \eqref{mi} and (1) that $\Hom _R (k, I)=0$.
Suppose that $I\not =0$.
Then there exists an $I$-regular element $(s, v)\in\m$ (cf. \cite[Proposition 1.2.3]{BH}).
It is easy to observe that the element $s\in\n$ is $S$-regular, contrary to the assumption that $\depth\, S=0$.
Therefore we have $I=0$, equivalently, $S=k$.
By (1) again, we obtain isomorphisms $R \cong k \ltimes k \cong k[[X]]/(X^2)$.
\qed
\end{pf}

The structure of an idealization of the form in the above proposition is uniquely determined if it has at least a non-free module of G-dimension zero.

\begin{cor}\label{corideal}
Let $(S, \n , k)$ be a local ring, $V \not = 0$ a finite-dimensional $k$-vector space, and $R = S \ltimes V$ the idealization of $V$ over $S$.
Then the following are equivalent:
\begin{enumerate}
\item[$(1)$]\ There is a non-free $R$-module in $\G (R)$;
\item[$(2)$]\ $R$ is Gorenstein;
\item[$(3)$]\ $R \cong k[[X]]/(X^2)$.
\end{enumerate}
\end{cor}

\begin{pf}
(3) $\Rightarrow$ (2): This implication is obvious.

(2) $\Rightarrow$ (1): Note that $\Kdim\, R = \depth\, R = \min\{\depth\, S, \depth _S V\} = 0$, namely, $R$ is an Artinian local ring.
Hence $k$ belongs to $\G (R)$.
Suppose that the $R$-module $k$ is free.
Then $R$ is regular, and hence $R$ is a field.
However, there is a non-zero element $v\in V$, and the element $(0, v)\in R$ is non-zero and nilpotent, which is a contradiction.
Thus $k$ is a non-free $R$-module in $\G (R)$.

(1) $\Rightarrow$ (2): Then, we see that there exists a non-free indecomposable $R$-module $M$ in $\G (R)$.
By definition it is reflexive.
Proposition \ref{ideal}(1) says that $M$ is isomorphic to $k$.
It follows that $R$ is Gorenstein.

(2) $\Rightarrow$ (3): Suppose that $\depth\, S>0$.
Then we especially have $\Kdim\, R = \Kdim\, S > 0$.
Since $\depth\, R = 0$, the local ring $R$ is not Cohen-Macaulay, and hence $R$ is not Gorenstein, which is a contradiction.
Therefore $\depth\, S=0$, and Proposition \ref{ideal}(2) implies that $R\cong k[[X]]/(X^2)$.
\qed
\end{pf}

%%%%% The first syzygy of the residue field (i.e. the maximal ideal) %%%%%
\subsection{The first syzygy of the residue field (i.e. the maximal ideal)}

The decomposability of the maximal ideal and the existence of a non-free module of G-dimension zero played essential roles in the achievement of Corollary \ref{corideal}.
From now on, we consider a local ring satisfying these conditions in more general settings.
First of all, let us describe the minimal free resolution of the residue class field of such a local ring.

\begin{prop}\label{efirst}
Let $(R, \m , k)$ be a local ring.
Suppose that there is a direct sum decomposition $\m = I \oplus J$ where $I, J$ are non-zero ideals of $R$.
Let $M$ be a non-free indecomposable module in $\G (R)$.
Then there exist $x, y\in\m$ such that
\begin{enumerate}
\item[$(1)$]\ $I=(x)$ and $J=(y)$,
\item[$(2)$]\ $(0:x)=(y)$ and $(0:y)=(x)$,
\item[$(3)$]\ $M$ is isomorphic to either $(x)$ or $(y)$.
\end{enumerate}
Hence the minimal free resolution of $k$ is as follows:
$$
\begin{CD}
\cdots @>{\left(
\begin{smallmatrix}
y & 0 \\
0 & x
\end{smallmatrix}
\right)}>> R^2 @>{\left(
\begin{smallmatrix}
x & 0 \\
0 & y
\end{smallmatrix}
\right)}>> R^2 @>{\left(
\begin{smallmatrix}
y & 0 \\
0 & x
\end{smallmatrix}
\right)}>> R^2 @>{\left(
\begin{smallmatrix}
x & y
\end{smallmatrix}
\right)}>> R @>>> k @>>> 0.
\end{CD}
$$
\end{prop}

\begin{pf}
Both $M^{\ast}$ and $\Omega M$ are also non-free indecomposable modules in $\G (R)$.
By virtue of Lemma \ref{nfi}, there are isomorphisms $M^{\ast}\cong\Hom _R (M, \m)=\Hom _R (M, I \oplus J)\cong\Hom _R (M, I)\oplus\Hom _R (M, J)$.
The indecomposability of $M^{\ast}$ implies that either $\Hom _R (M, I)=0$ or $\Hom _R (M, J)=0$.
We may assume that
\begin{equation}\label{assume}
\Hom _R (M, J)=0.
\end{equation}

There is an exact sequence
\begin{equation}\label{syz}
0 \to \Omega M \to R^n \to M \to 0.
\end{equation}
Dualizing this by $J$, we obtain another exact sequence $\Hom _R (M, J) \to J^n \to \Hom _R (\Omega M, J)$.
We have $\Hom _R (\Omega M, J)\not =0$ by \eqref{assume}.
Applying the above argument to the module $\Omega M$ yields
\begin{equation}\label{ni}
\Hom _R (\Omega M, I)=0.
\end{equation}
Also, dualizing \eqref{syz} by $I$, we get an exact sequence $0 \to \Hom _R (M, I) \to I^n \to \Hom _R (\Omega M, I)$, and hence $M^{\ast}\cong\Hom _R (M, I)\cong I^n$.
The indecomposability of $M^{\ast}$ implies that $n=1$ (i.e. $M$ is cyclic), and $M^{\ast}\cong I$.
Let $\alpha : M^{\ast}\to I$ denote this isomorphism, and write $M=Rz$ for some $z\in M$.
Then it is easy to check that $\alpha$ is a map defined by $\alpha (\sigma )=\sigma (z)$ for $\sigma\in M^{\ast}$.

We also have $M \cong M^{\ast\ast}\cong \Hom _R (M^{\ast}, \m ) \cong \Hom _R (M^{\ast}, I)\oplus\Hom _R (M^{\ast}, J)$.
Note that $\Hom _R (M^{\ast}, I)$ is isomorphic to $\Hom _R (I, I)$, which contains the identity map of $I$.
Hence $\Hom _R (M^{\ast}, I)\not =0$ and therefore $\Hom _R (M^{\ast}, J)=0$.
Applying the above argument to the module $M^{\ast}$, we see that $M^{\ast}$ is also cyclic and $M\cong M^{\ast\ast}\cong I$.
Thus, we have shown that $M\cong M^{\ast}\cong I$ and these modules are cyclic.
Noting \eqref{ni} and applying the above argument to the module $\Omega M$, we see that $\Omega M\cong (\Omega M)^{\ast}\cong J$ and these modules are cyclic.

Write $I=(x)$ and $J=(y)$.
Then $M$ is isomorphic to the principal ideal $(x)$.
Apply the above argument to $(x)$ instead of $M$, and we have an isomorphism $\alpha : (x)^{\ast}\to (x)$ which is defined by $\alpha (\sigma )=\sigma (x)$ for $\sigma\in (x)^{\ast}$.
Consider a composite map $(0:(0:x)) \overset{\gamma}{\to} (R/(0:x))^{\ast} \overset{\beta}{\to} (x)^{\ast} \overset{\alpha}{\to} (x)$, where $\beta , \gamma$ are natural isomorphisms.
We easily see that this composite map is the identity map.
Hence $(0:(0:x))=(x)$.
Similarly, we also have $(0:(0:y))=(y)$.
Since $(0:x) = \Omega (x) \cong \Omega M \cong (y)$, we have $(x) = (0:(0:x))=\Ann _R (0:x)= \Ann _R (y)=(0:y)$, and therefore $(0:x)=\Ann _R (x)=\Ann _R (0:y) =(0:(0:y))= (y)$.
Thus we obtain the minimal free resolutions of $(x)$ and $(y)$:
$$
\begin{cases}
\cdots\overset{y}{\to} R \overset{x}{\to} R \overset{y}{\to} R \to (x) \to 0,\\
\cdots\overset{x}{\to} R \overset{y}{\to} R \overset{x}{\to} R \to (y) \to 0.
\end{cases}
$$
Taking the direct sum of these exact sequence, we get
$$
\begin{CD}
\cdots @>{\left(
\begin{smallmatrix}
y & 0 \\
0 & x
\end{smallmatrix}
\right)}>> R^2 @>{\left(
\begin{smallmatrix}
x & 0 \\
0 & y
\end{smallmatrix}
\right)}>> R^2 @>{\left(
\begin{smallmatrix}
y & 0 \\
0 & x
\end{smallmatrix}
\right)}>> R^2 @>>> \m @>>> 0.
\end{CD}
$$
Joining this to the natural exact sequence $0 \to \m \to R \to k \to 0$ constructs the minimal free resolution of $k$ in the assertion.
\qed
\end{pf}

We denote by $\edim R$ the embedding dimension of a local ring $R$.
When a homomorphic image of a regular local ring is given, we can choose a minimal presentation of the ring in the following sense:

\begin{prop}\label{pres}
Let $R$ be a homomorphic image of a regular local ring.
Then there exist a regular local ring $(S, \n )$ and an ideal $I$ of $S$ contained in $\n ^2$ such that $R\cong S/I$.
\end{prop}

Here we introduce a famous result due to Tate \cite[Theorem 6]{Tate}.
See also \cite[Remarks 8.1.1(3)]{Avramov}.

\begin{lem}[Tate]\label{tate}
Let $(S, \n , k)$ be a regular local ring, $I$ an ideal of $S$ contained in $\n ^2$, and $R=S/I$ a residue class ring.
Suppose that the complexity of $k$ over $R$ is at most one.
(In other words, the set of all the Betti numbers of the $R$-module $k$ is bounded.)
Then $I$ is a principal ideal.
\end{lem}

We denote by $\beta _i ^R (M)$ the $i$th Betti number of a module $M$ over a local ring $R$.
Handling the above results, we can determine the structure of a local ring with decomposable maximal ideal having a non-free module of G-dimension zero, as follows:

\begin{thm}\label{first}
Let $(S, \n, k)$ be a regular local ring, $I$ an ideal of $S$ contained in $\n ^2$, and $R=S/I$ a residue class ring.
Suppose that there exists a non-free $R$-module in $\G (R)$.
Then the following conditions are equivalent:
\begin{enumerate}
\item[$(1)$]\ The maximal ideal of $R$ is decomposable;
\item[$(2)$]\ $\Kdim\, S=2$ and $I=(xy)$ for some regular system of parameter $x, y$ of $S$.
\end{enumerate}
\end{thm}

\begin{pf}
Let $\m = \n /I$ be the maximal ideal of $R$.

(2) $\Rightarrow$ (1): It is easy to see that $\m = xR\oplus yR$ and that $xR, yR$ are non-zero.

(1) $\Rightarrow$ (2): First of all, note from the condition (1) that $R$ is not an integral domain, hence is not a regular local ring.
Proposition \ref{efirst} says that $\m = xR\oplus yR$ for some $x, y\in\n$, and that $\beta _i ^R (k)=2$ for every $i\geq 2$.
It follows from Lemma \ref{tate} that $I$ is a principal ideal.
Hence $R$ is a hypersurface.
We write $I=(f)$ for some $f\in\n ^2$.
Since $\m$ is decomposable, the local ring $R$ is not Artinian.
(Over an Artinian Gorenstein local ring, the intersection of non-zero ideals is also non-zero; cf. \cite[Exercise 3.2.15]{BH}.)
Hence we have $0<\Kdim\, R<\edim R=2$, which says that $\Kdim\, R=1$ and $\Kdim\, S=2$.

Note that $\n = (x, y, f)$.
Because $\edim S=\Kdim\, S=2$, one of the elements $x, y, f$ belongs to the ideal generated by the other two elements.
Noting that the images of elements $x, y$ in $\m$ form a minimal system of generators of $\m$, we see that $f\in (x, y)$, and hence $x, y$ is a regular system of parameters of $S$.
On the other hand, noting $xR\cap yR=0$, we get $xy\in I=(f)$.
Write $xy=cf$ for some $c\in S$.
Since the associated graded ring $\gr _{\n} (S)$ is a polynomial ring over $k$ in two variables $\overline{x}, \overline{y}\in\n /{\n}^2$, we especially have $\overline{x}\overline{y}\not =0$ in ${\n}^2 /{\n}^3$, namely, $xy\not\in\n ^3$.
It follows that $c\not\in\n$ because $f\in\n ^2$.
Therefore the element $c$ is a unit of $S$, and thus $I=(xy)$.
\qed
\end{pf}

Using Theorem \ref{first} and Cohen's structure theorem, we obtain the following corollary.

\begin{cor}\label{corfirst}
Let $(R, \m )$ be a complete local ring.
The following conditions are equivalent:
\begin{enumerate}
\item[$(1)$]\ There is a non-free module in $\G (R)$, and $\m$ is decomposable;
\item[$(2)$]\ $R$ is Gorenstein, and $\m$ is decomposable;
\item[$(3)$]\ There are a complete regular local ring $S$ of dimension two and a regular system of parameters $x, y$ of $S$ such that $R\cong S/(xy)$.
\end{enumerate}
\end{cor}

Note that the finiteness of G-dimension is independent of completion.
Thus, Corollary \ref{corfirst} not only gives birth to a generalization of \cite[Proposition 2.3]{Takahashi} but also guarantees that Question \ref{ques} is true if $n=1$.

%%%%% The second syzygy of the residue field %%%%%
\subsection{The second syzygy of the residue field}

As far as here, we have observed a local ring whose maximal ideal is decomposable.
From here to the end of this paper, we will observe a local ring such that the second syzygy module of the residue class field is decomposable.
We begin with the following theorem, which implies that Question \ref{ques} is true if $n=2$.

\begin{thm}\label{second}
Let $(R, \m, k)$ be a local ring.
Suppose that $\m$ is indecomposable and that $\Omega _R ^2 k$ has a non-zero proper direct summand of finite G-dimension.
Then $R$ is a Gorenstein ring of dimension two.
\end{thm}

\begin{pf}
Replacing $R$ with its $\m$-adic completion, we may assume that $R$ is a complete local ring.
In particular, note that $R$ is Henselian.
We have $\Omega _R ^2 k=M\oplus N$ for some non-zero $R$-modules $M$ and $N$ with $\g _R M<\infty$.
There is an exact sequence $0 \longrightarrow M\oplus N \overset{(f, g)}{\longrightarrow} R^e \longrightarrow \m \longrightarrow 0$ of $R$-modules, where $e=\edim R$.
Setting $A=\Coker f$ and $B=\Coker g$, we get exact sequences
\begin{equation}\label{coker}
\begin{cases}
0 \to M \overset{f}{\to} R^e \overset{\alpha}{\to} A \to 0,\\
0 \to N \overset{g}{\to} R^e \overset{\beta}{\to} B \to 0.
\end{cases}
\end{equation}
It is easily observed that there are exact sequences
\begin{equation}\label{key}
0 \longrightarrow R^e \overset{\binom{\alpha}{\beta}}{\longrightarrow} A\oplus B \longrightarrow \m \longrightarrow 0
\end{equation}
and
\begin{equation}\label{snake}
\begin{cases}
0 \longrightarrow M \overset{\beta f}{\longrightarrow} B \longrightarrow \m \longrightarrow 0,\\
0 \longrightarrow N \overset{\alpha g}{\longrightarrow} A \longrightarrow \m \longrightarrow 0.
\end{cases}
\end{equation}

\begin{claim}\label{ext2}
We have $\Ext _R ^2 (k, R)\not = 0$.
(Hence $\depth\, R\leq 2$.)
\end{claim}

\begin{pf}
Suppose that $\Ext _R ^2 (k, R)=0$.
Then $\Ext _R ^1 (\m , R^e)\cong\Ext _R ^2 (k, R^e)=0$.
Hence the exact sequence \eqref{key} splits, and therefore we have an isomorphism $A\oplus B\cong R^e\oplus\m$.
Since the maximal ideal $\m$ is indecomposable, it follows from the Krull-Schmidt theorem that $\m$ is isomorphic to a direct summand of $A$ or $B$.
If $\m$ is isomorphic to a direct summand of $A$, then $B$ is isomorphic to a direct summand of $R^e$.
Hence $B$ is a free $R$-module of rank at most $e$.
Denote by $b$ the rank of $B$.
Since the second sequence in \eqref{coker} splits, the $R$-module $N$ is a free module of rank $e-b$.
Noting that there is a surjective homomorphism from $B$ to $\m$ by \eqref{snake}, we have $b=\nu _R (B)\geq\nu _R (\m)=e$.
This means that $b=e$, and hence $N=0$, which is a contradiction.
We can get a contradiction along the same lines in the case where $\m$ is isomorphic to a direct summand of $B$.
Thus, we obtain $\Ext _R ^2 (k, R)\not =0$.
\qed
\end{pf}

Fix a non-free indecomposable module $X\in\G (R)$.
Applying the functor $\Hom _R (X, -)$ to \eqref{key} gives an exact sequence $0 \to (X^{\ast})^e \to \Hom _R (X, A)\oplus\Hom _R (X, B) \to \Hom _R (X, \m ) \to 0$ and an isomorphism
\begin{equation}\label{isom}
\Ext _R ^1 (X, A)\oplus\Ext _R ^1 (X, B)\cong\Ext _R ^1 (X, \m).
\end{equation}
We have $(X^{\ast})^e\in\G (R)$ and $\Hom _R (X, \m )\in\G (R)$ by Lemma \ref{nfi}, hence $\Hom _R (X, A)\in\G (R)$.

Take the first syzygy module of $X$; we have an exact sequence $0 \to \Omega X \to R^n \to X \to 0$.
Dualizing this sequence by $A$, we obtain an exact sequence $0 \to \Hom _R (X, A) \to A^n \to \Hom _R (\Omega X, A) \to \Ext _R ^1 (X, A) \to 0$.
Divide this into two short exact sequences
\begin{equation}\label{divide}
\begin{cases}
0 \to \Hom _R (X, A) \to A^n \to C \to 0, \\
0 \to C \to \Hom _R (\Omega X, A) \to \Ext _R ^1 (X, A) \to 0
\end{cases}
\end{equation}
of $R$-modules.
Since $\Omega X$ is also a non-free indecomposable module in $\G (R)$, applying the above argument to $\Omega X$ instead of $X$ shows that the module $\Hom _R (\Omega X, A)$ also belongs to $\G (R)$.
We have $\g _R (A^n)<\infty$ by the first sequence in \eqref{coker}.
Hence it follows from \eqref{divide} that $\g _R C<\infty$, and
\begin{equation}\label{gdimext}
\g _R (\Ext _R ^1 (X, A)) < \infty.
\end{equation}

On the other hand, applying the functor $\Hom _R (X, -)$ to the natural exact sequence $0 \to \m \to R \to k \to 0$, we get an exact sequence $0 \to \Hom _R (X, \m ) \to X^{\ast} \to \Hom _R (X, k) \to \Ext _R ^1 (X, \m ) \to 0$.
Lemma \ref{nfi} implies that $\Hom _R (X, k)\cong\Ext _R ^1 (X, \m )$, hence $\Ext _R ^1 (X, \m )$ is a $k$-vector space.
Since $\Ext _R ^1 (X, A)$ is contained in $\Ext _R ^1 (X, \m )$ by \eqref{isom}, 
\begin{equation}\label{extvec}
\Ext _R ^1 (X, A)\text{ is a }k\text{-vector space}.
\end{equation}

\begin{claim}\label{goren}
The local ring $R$ is Gorenstein.
\end{claim}

\begin{pf}
Suppose that $R$ is not Gorenstein.
Then we must have $\Ext _R ^1 (G, A)=0$ for any $G\in\G (R)$ by \eqref{gdimext} and \eqref{extvec}.
We have an exact sequence
\begin{equation}\label{cosyzygyx}
0 \to X \to R^m \to \Omega ^{-1} X \to 0,
\end{equation}
and note that $\Omega ^{-1} X$ belongs to $\G (R)$.
The exact sequences \eqref{cosyzygyx} and \eqref{coker} yield isomorphisms $\Ext _R ^1 (X, M) \cong \Ext _R ^2 (\Omega ^{-1} X, M)\cong \Ext _R ^1 (\Omega ^{-1} X, A)= 0$.
This means that
\begin{equation}\label{extm}
\Ext _R ^1 (G, M)=0
\end{equation}
for any $G\in\G (R)$.
On the other hand, since $\depth _R M \geq \depth _R (\Omega ^2 k) \geq \min\{ 2, \depth\, R\}= \depth\, R$ by \cite[Exercise 1.3.7]{BH} and Claim \ref{ext2}, $M$ belongs to $\G (R)$.
Hence there is an exact sequence of the form $0 \to M \to R^l \to \Omega ^{-1} M \to 0$, and this splits because $\Ext _R ^1 (\Omega ^{-1} M, M)=0$ by \eqref{extm}.
Thus the $R$-module $M$ is free.
Theorem \ref{dutta} implies that $R$ is regular, which contradicts our assumption that $R$ is not Gorenstein.
This contradiction proves the claim.
\qed
\end{pf}

Since the only number $i$ such that $\Ext _R ^i (k, R)\not =0$ is the Krull dimension of $R$ if $R$ is Gorenstein, it follows from the above two claims that $R$ is a Gorenstein local ring of dimension two, which completes the proof of the theorem.
\qed
\end{pf}

The above theorem interests us in the investigation of a Gorenstein local ring of dimension two such that the second syzygy module of the residue class field is decomposable.
Our result concerning this is stated as follows.

\begin{thm}\label{2dimgor}
Let $(S, \n , k)$ be a regular local ring, $I$ an ideal of $S$ contained in $\n ^2$, and $R=S/I$ a residue class ring.
Suppose that $R$ is a Henselian Gorenstein ring of dimension two.
Then the following are equivalent:
\begin{enumerate}
\item[$(1)$]\ $\Omega _R ^2 k$ is decomposable;
\item[$(2)$]\ $\Kdim\, S=3$ and $I=(xy-zf)$ for some regular system of parameters $x, y, z$ of $S$ and $f\in\n$.
\end{enumerate}
\end{thm}

It is necessary to prepare three elementary lemmas to prove this theorem.
The first and third ones are both well-known and easy to check, and we omit the proofs.

\begin{lem}\label{hs}
Let $(S, \n , k)$ be a regular local ring of dimension three and $R=S/(f)$ a hypersurface with $f\in\n ^2$.
Then $f=x f_x + y f_y + z f_z$ for some $f_x , f_y , f_z \in\n$, and the minimal free resolution of $k$ over $R$ is as follows:
$$
\cdots\overset{C}{\longrightarrow} R^4 \overset{D}{\longrightarrow} R^4 \overset{C}{\longrightarrow} R^4 \overset{D}{\longrightarrow} R^4 \overset{C}{\longrightarrow} R^4 \overset{B}{\longrightarrow} R^3 \overset{A}{\longrightarrow} R \longrightarrow k \longrightarrow 0,
$$
where
$$
\begin{array}{rccccl}
A & = & \left(
\begin{smallmatrix}
 x & y & z
\end{smallmatrix}
\right), & B & = & \left(
\begin{smallmatrix}
 0 & -z &  y & f_x \\
 z &  0 & -x & f_y \\
-y &  x &  0 & f_z
\end{smallmatrix}
\right),\\
\\
C & = & \left(
\begin{smallmatrix}
   0 & -f_z &  f_y & x \\
 f_z &    0 & -f_x & y \\
-f_y &  f_x &    0 & z \\
  -x &   -y &   -z & 0
\end{smallmatrix}
\right), & D & = & \left(
\begin{smallmatrix}
   0 &   -z &    y & f_x \\
   z &    0 &   -x & f_y \\
  -y &    x &    0 & f_z \\
-f_x & -f_y & -f_z &   0
\end{smallmatrix}
\right).
\end{array}
$$
\end{lem}

\begin{lem}\label{wl}
Let $(R, \m, k)$ be a local ring and $x\in\m - \m ^2$ an $R$-regular element.
Then we have a split exact sequence $0 \to k \overset{\theta}{\to} \m /x\m \overset{\pi}{\to} \m /xR \to 0$, where $\theta$ is defined by $\theta (\overline{a})=\overline{xa}$ for $\overline{a}\in R/\m =k$ and $\pi$ is the natural surjection.
\end{lem}

\begin{pf}
Let $x_1, x_2, \dots, x_n$ be a minimal system of generators of $\m$ with $x_1=x$.
Define a homomorphism $\varepsilon : \m /x\m \to k$ by $\varepsilon (\sum _{i=1} ^n \overline{x_i a_i})=\overline{a_1}$.
We easily see that the composite map $\varepsilon\theta$ is the identity map of $k$, which means that $\theta$ is a split-monomorphism.
\qed
\end{pf}

\begin{lem}\label{dvr}
Let $(R, \m, k)$ be a Cohen-Macaulay local ring of dimension one.
Then $R$ is a discrete valuation ring if and only if $\m ^{\ast}$ is a cyclic $R$-module.
\end{lem}

Now let us prove Theorem \ref{2dimgor}.

\begin{tpf}
(2) $\Rightarrow$ (1): We have $xy-zf=x\cdot 0+y\cdot x+z\cdot (-f)$.
Lemma \ref{hs} gives a finite free presentation $R^4 \overset{C}{\longrightarrow} R^4 \longrightarrow \Omega _R ^2 k \longrightarrow 0$ of the $R$-module $\Omega _R ^2 k$, where $C=\left(
\begin{smallmatrix}
 0 &  f &  x &  x \\
-f &  0 &  0 &  y \\
-x &  0 &  0 &  z \\
-x & -y & -z &  0
\end{smallmatrix}
\right)$.
Putting $P=\left(
\begin{smallmatrix}
1 & 0 & 0 &  0 \\
0 & 0 & 1 & -1 \\
0 & 0 & 1 &  0 \\
0 & 1 & 0 &  0
\end{smallmatrix}
\right)$ and $Q=\left(
\begin{smallmatrix}
0 & 0 & -1 &  0 \\
0 & 1 &  0 &  0 \\
1 & 0 &  0 & -1 \\
0 & 0 &  0 &  1
\end{smallmatrix}
\right)$, we obtain $PCQ=
\begin{pmatrix}
U & 0\\
0 & {}^t U
\end{pmatrix}
$, where $U=\left(
\begin{smallmatrix}
x & f\\
z & y
\end{smallmatrix}
\right)$.
It is easily seen that the matrices $P, Q$ are invertible.
Denoting by $M$ (resp. $N$) the cokernel of the homomorphism defined by the matrix $U$ (resp. ${}^t U$), we get an isomorphism $\Omega _R ^2 k \cong M\oplus N$.

(1) $\Rightarrow$ (2): First of all, note that the local ring $R$ is not regular.
We denote by $\m$ the maximal ideal $\n/I$ of $R$.

Suppose that there exists an element $z\in\n -\n ^2$ whose image in $\m$ is an $R$-regular element such that the module $\m / zR$ is decomposable.
Then the assertion (2) follows.
Indeed, put $\overline{(-)}=(-)\otimes _S S/(z)$.
Note that $\overline{S}$ is also a regular local ring because $z$ is a minimal generator of the maximal ideal $\n$ of $S$ (see the proof of Proposition \ref{pres}).
Since the maximal ideal $\m\overline{R}$ of $\overline{R}$ is decomposable, we can apply Theorem \ref{first} and see that $\Kdim\,\overline{S}=2$ and $I\overline{S}=xy\overline{S}$ for some $x, y\in\n$ whose images in $\overline{S}$ form a regular system of parameter of $\overline{S}$.
Hence $\overline{R}=\overline{S}/xy\overline{S}$ is a hypersurface, in particular a complete intersection, of dimension one.
Therefore $R$ is a complete intersection of dimension two by \cite[Theorem 2.3.4(a)]{BH}.
Since $S$ is a regular local ring of dimension three with regular system of parameter $x, y, z$, the ideal $I$ is generated by an $S$-sequence by \cite[Theorem 2.3.3(c)]{BH}.
Noting $\height I = \Kdim\, S - \Kdim\, R = 1$, we see that $I$ is a principal ideal.
Write $I=(l)$ for some $l\in I$.
There is an element $f\in S$ such that $l=xy-zf$.
Assume that $f\not\in\n$.
Then $f$ is a unit of $S$, and we see that $zR\subseteq xyR$.
Hence $\m = (x, y)R$, and $\edim R=\Kdim\, R=2$.
This implies that $R$ is regular, which is a contradiction.
It follows that $f\in\n$.

On the other hand, if $z\in\n$ is an element whose image in $\m$ is $R$-regular such that $\m /zR$ is decomposable, then $z\not\in\n ^2$.
Indeed, assume $z\in\n ^2$.
Then we have $I+(z)\subseteq\n ^2$.
Since $R/zR=S/I+(z)$, it follows from Theorem \ref{first} that $\Kdim\,S=2$.
Since $\Kdim\,R=2$, we have $I=0$, equivalently $R=S$.
In particular $R$ is regular, which is a contradiction.

Thus, it suffices to show the existence of an $R$-regular element $w\in\m$ such that $\m / (w)$ is decomposable.
Let $E$ denote the fundamental module of $R$.
Proposition \ref{fund}(4) says that we can write $E=M\oplus N$ for some non-zero $R$-modules $M$ and $N$.
Hence the fundamental sequence of $R$ is as follows:
\begin{equation}\tag{a}\label{fundseq}
0 \longrightarrow R \overset{\binom{\sigma}{\tau}}{\longrightarrow} M\oplus N \overset{(f, g)}{\longrightarrow} \m \longrightarrow 0.
\end{equation}

Take an $R$-regular element $w\in\m - \m ^2$, and set $\overline{(-)}=(-)\otimes _R R/(w)$.
If $\m\overline{R}$ is decomposable, then our aim is attained.
Hence let $\m\overline{R}$ be indecomposable.
The sequence \eqref{fundseq} induces another exact sequence $0 \longrightarrow \overline{R} \overset{\binom{\overline{\sigma}}{\overline{\tau}}}{\longrightarrow} \overline{M}\oplus\overline{N} \overset{(\overline{f}, \overline{g})}{\longrightarrow} \overline{\m} \longrightarrow 0$.
(Here, the injectivity of the map $\binom{\overline{\sigma}}{\overline{\tau}}$ follows from the fact that $w$ is an $\m$-regular element.)
According to Lemma \ref{wl}, the natural surjection $\pi : \overline{\m} \to \m\overline{R}$ is a split-epimorphism with kernel isomorphic to $k$.
Hence there exists a split-monomorphism $\rho : \m\overline{R} \to \overline{\m}$ such that $\pi\rho =1$.
Then note that the cokernel of $\rho$ is isomorphic to $k$.
On the other hand (cf. Proposition \ref{fund}), the homomorphism $(\overline{f}, \overline{g})$ is a $\G (\overline{R})$-precover of $\overline{\m}$.
Therefore there exists a homomorphism $\binom{\alpha}{\beta} : \m\overline{R} \to \overline{M}\oplus\overline{N}$ such that $\rho = (\overline{f}, \overline{g})\binom{\alpha}{\beta} = \overline{f}\alpha + \overline{g}\beta$.
Set $e=\edim R$, $m=\nu _R (M)$, and $n=\nu _R (N)$.

\setcounter{claim}{0}
\begin{claim}\label{mn}
We have either $(m,n)=(e-1,2)$ or $(m,n)=(2,e-1)$.
\end{claim}

\begin{pf}
Since $\rho$ is a split-monomorphism, so is the homomorphism $\binom{\alpha}{\beta}$.
There is a commutative diagram
$$
\begin{CD}
@. @. 0 @. 0 \\
@. @. @VVV @VVV \\
@. @. \m\overline{R} @= \m\overline{R} \\
@. @. @V{\left(
\begin{smallmatrix}
\alpha\\
\beta
\end{smallmatrix}
\right)}VV @V{\rho}VV \\
0 @>>> \overline{R} @>{\left(
\begin{smallmatrix}
\overline{\sigma}\\
\overline{\tau}
\end{smallmatrix}
\right)}>> \overline{M}\oplus\overline{N} @>{\left(
\begin{smallmatrix}
\overline{f} & \overline{g}
\end{smallmatrix}
\right)}>> \overline{\m} @>>> 0 \\
@. @| @VVV @VVV \\
0 @>>> \overline{R} @>>> C @>>> k @>>> 0 \\
@. @. @VVV @VVV \\
@. @. 0 @. 0
\end{CD}
$$
of $\overline{R}$-modules with exact rows and columns, and we have an isomorphism $\overline{M}\oplus\overline{N}\cong\m\overline{R}\oplus C$.
The indecomposability of $\m\overline{R}$ and the Krull-Schmidt theorem yield that $\m\overline{R}$ is isomorphic to a direct summand of either $\overline{M}$ or $\overline{N}$.

Let us consider the case where $\m\overline{R}$ is isomorphic to a direct summand of $\overline{M}$.
There is an $\overline{R}$-module $L$ such that $\overline{M}\cong\m\overline{R}\oplus L$.
The Krull-Schmidt theorem again yields an isomorphism
\begin{equation}\tag{b}\label{c1}
C\cong\overline{N}\oplus L.
\end{equation}
Note that $\overline{N}$ and $L$ are isomorphic to direct summands of $\overline{E}$.
Proposition \ref{fund} implies that the $\overline{R}$-module $\overline{E}$ belongs to $\G (\overline{R})$.
The $\overline{R}$-modules $\overline{N}, L$ also belong to $\G (\overline{R})$, and so does $C$ by \eqref{c1}.
Therefore the exact sequence
\begin{equation}\tag{c}\label{no1}
0 \to \overline{R} \to C \to k \to 0
\end{equation}
in the above diagram does not split because $\depth\,C=1>0$.
On the other hand, noting that $\overline{R}$ is a Gorenstein local ring of dimension one, we have $\Hom _{\overline{R}} (k, \overline{R})=0$ and $\Ext _{\overline{R}} ^1 (k, \overline{R})\cong k$.
Dualizing the natural exact sequence $0 \to \m\overline{R} \to \overline{R} \to k \to 0$, we have another exact sequence
\begin{equation}\tag{d}\label{no2}
0 \to \overline{R} \to \Hom _{\overline{R}} (\m\overline{R}, \overline{R})  \to k \to 0.
\end{equation}
Note that the maximal ideal $\m\overline{R}$ of $\overline{R}$ belongs to $\G (\overline{R})$, hence so does $\Hom _{\overline{R}} (\m\overline{R}, \overline{R})$.
Therefore the exact sequence \eqref{no2} does not split because $\depth\,\Hom _{\overline{R}} (\m\overline{R}, \overline{R})=1>0$.

Thus, we have obtained two non-split exact sequences \eqref{no1} and \eqref{no2} of $\overline{R}$-modules.
Since $\Ext _{\overline{R}} ^1 (k, \overline{R})\cong k$, we obtain an isomorphism
\begin{equation}\tag{e}\label{c2}
C\cong\Hom _{\overline{R}} (\m\overline{R}, \overline{R}).
\end{equation}
The isomorphisms \eqref{c1} and \eqref{c2} give other isomorphisms $\m\overline{R} \cong \Hom _{\overline{R}} (\Hom _{\overline{R}} (\m\overline{R}, \overline{R}), \overline{R})\cong \Hom _{\overline{R}} (\overline{N}\oplus L, \overline{R})\cong \Hom _{\overline{R}} (\overline{N}, \overline{R})\oplus\Hom _{\overline{R}} (L, \overline{R})$.
Note that $\overline{N}\not =0$ and $L$ are reflexive $\overline{R}$-modules, hence $\Hom _{\overline{R}} (\overline{N}, \overline{R})\not =0$.
Since $\m\overline{R}$ is indecomposable, we have $\Hom _{\overline{R}} (L, \overline{R})=0$, and hence $L=0$.
Thus we get two isomorphisms $\overline{M}\cong\m\overline{R}$ and $\overline{N}\cong\Hom _{\overline{R}} (\m\overline{R}, \overline{R})$.
Therefore $m=\nu _{\overline{R}} (\overline{M})=\edim\overline{R}=e-1$ because $w\not\in\m ^2$, and $n=\nu _{\overline{R}} (\overline{N})=\nu _{\overline{R}} (\Hom _{\overline{R}} (\m\overline{R}, \overline{R}))$.
Lemma \ref{dvr} implies that $n\geq 2$.
On the other hand, it follows from the fundamental sequence \eqref{fundseq} that $m+n= \nu _R (M\oplus N)\leq \nu _R (R)+\nu _R (\m ) = 1+e$.
Hence we see that $n=2$.

In the case where $\m\overline{R}$ is isomorphic to a direct summand of $\overline{N}$, a similar argument yields $m=2$ and $n=e-1$.
\qed
\end{pf}

On the other hand, we have $1 = \pi\rho = \pi\overline{f}\alpha + \pi\overline{g}\beta$ in $\End _{\overline{R}} (\m\overline{R})$.
Since $\m\overline{R}$ is indecomposable, the endomorphism ring $\End _{\overline{R}} (\m\overline{R})$ is a local ring (cf. \cite[Proposition (1.18)]{Yoshino}), and hence either $\pi\overline{f}\alpha$ or $\pi\overline{g}\beta$ is a unit of this ring, in other words, is an automorphism.
Put $\a = \Im f$ and $\b = \Im g$.

\begin{claim}\label{2kinds}
If $\pi\overline{f}\alpha$ (resp. $\pi\overline{g}\beta$) is an automorphism, then $\m = \a + (w)$ (resp. $\m = \b + (w)$) and $\grade\,\a >0$ (resp. $\grade\,\b >0$).
\end{claim}

\begin{pf}
Suppose that $\pi\overline{f}\alpha$ is an automorphism.
Then $\pi\overline{f}$ is a split-epimorphism, and so in particular a surjection.
Hence $\m\overline{R} = \a\overline{R}$, and therefore $\m = \a + (w)$.
There exists an $\overline{R}$-regular element in $\m\overline{R}=\a\overline{R}$.
We can choose an element $v\in\a$ whose image in $\m\overline{R}$ is $\overline{R}$-regular.
Since $w, v$ is an $R$-regular sequence, so is the sequence $v, w$.
Thus $v$ is an $R$-regular element.
The proof of the other case is similar.
\qed
\end{pf}

\begin{claim}\label{grade}
We have both $\grade\,\a >0$ and $\grade\,\b >0$.
\end{claim}

\begin{pf}
It is enough to show the claim only in the case where $\pi\overline{f}\alpha$ is an automorphism.
Then Claim \ref{2kinds} says that $\m = \a + (w)$ and $\grade\,\a >0$.
Take an $R$-regular element $v\in\a - \m ^2$.
Applying the above argument to the element $v$ instead of $w$, we see that either of the following holds:
\begin{enumerate}
\item[(i)]\ $\m = \a + (v)$ and $\grade\,\a >0$;
\item[(ii)]\ $\m = \b + (v)$ and $\grade\,\b >0$.
\end{enumerate}
However, if the statement (i) holds, then we have $\m = \a$, which means that the homomorphism $f : M \to \m$ is surjective.
Hence $m=\nu _R (M)\geq\nu _R (\m )=e$.
It follows from Claim \ref{mn} that $m=2$, and hence $e\leq 2$.
But this can not happen because $R$ is a non-regular local ring of dimension two.
Consequently the statement (ii) must hold, and we obtain $\grade\,\b >0$, as desired.
\qed
\end{pf}

Put $x=\sigma (1)$ and $y=\tau (1)$.
Then $f(x)+g(y)=(f, g)\binom{\sigma}{\tau} (1)=0$.
Set $v=f(x)=-g(y)\in\a\cap\b$.
Take an element $a\in\a\cap\b$.
Then we have $a=f(p)=g(q)$ for some $p\in M$ and $q\in N$.
Hence $\binom{p}{-q}\in\Ker (f,g)=\Im\binom{\sigma}{\tau}$, and therefore $\binom{p}{-q} = b\binom{x}{y}$ for some $b\in R$.
Thus $p=bx$, and we get $a=f(p)=f(bx)=bv\in (v)$.
It follows that $\a\cap\b = (v)$.
Since $\grade (v) = \grade (\a\cap\b ) = \inf\{\grade\,\a , \grade\,\b\} > 0$ by \cite[Proposition 1.2.10(c)]{BH} and Claim \ref{grade}, the element $v$ is an $R$-regular element.

Set $\overline{(-)}=(-)\otimes _R R/(v)$.
Since $\a + \b = \m$ and $\a \cap \b = (v)$, there is a natural exact sequence $\omega : 0 \to \overline{R} \to R/\a \oplus R/\b \to k \to 0$ of $\overline{R}$-modules.
Suppose that this exact sequence splits.
Then we have an isomorphism $R/\a \oplus R/\b \cong \overline{R} \oplus k$, and it is seen from the Krull-Schmidt theorem that $k$ is isomorphic to either $R/\a$ or $R/\b$.
Hence we have either $\m = \a$ or $\m = \b$, and the same argument as the end of the proof of Claim \ref{grade} yields a contradiction.
Thus the exact sequence $\omega$ does not split.

On the other hand, dualizing the natural exact sequence $0 \to \m\overline{R} \to \overline{R} \to k \to 0$, we have a non-split exact sequence $0 \to \overline{R} \to \Hom _{\overline{R}} (\m\overline{R}, \overline{R}) \to k \to 0$.
Since $\Ext _{\overline{R}} ^1 (k, \overline{R})\cong k$, we obtain an isomorphism $R/\a \oplus R/\b \cong \Hom _{\overline{R}} (\m\overline{R}, \overline{R})$, and $\Hom _{\overline{R}} (\m\overline{R}, \overline{R})$ belongs to $\G (\overline{R})$.
It follows that both $R/\a$ and $R/\b$ belong to $\G (\overline{R})$, hence they are reflexive over $\overline{R}$.
Therefore the $\overline{R}$-dual modules $\Hom _{\overline{R}} (R/\a , \overline{R})$ and $\Hom _{\overline{R}} (R/\b , \overline{R})$ are non-zero, which proves that $\m\overline{R}$ is decomposable.
This completes the proof of our theorem.
\qed
\end{tpf}

Combining Theorem \ref{second} with Theorem \ref{2dimgor} gives birth to the following corollary.
Compare it with Corollary \ref{corfirst}.

\begin{cor}
Let $(R, \m , k)$ be a complete local ring.
Suppose that $\m$ is indecomposable.
Then the following conditions are equivalent:
\begin{enumerate}
\item[$(1)$]\ $\Omega _R ^2 k$ has a non-zero proper direct summand of finite G-dimension;
\item[$(2)$]\ $R$ is Gorenstein, and $\Omega _R ^2 k$ is decomposable;
\item[$(3)$]\ There are a complete regular local ring $(S, \n )$ of dimension three, a regular system of parameters $x, y, z$ of $S$, and $f\in\n$ such that $R\cong S/(xy-zf)$.
\end{enumerate}
\end{cor}

Lastly, we recall a result of Yoshino and Kawamoto, which is related to Theorem \ref{2dimgor}.
A homomorphic image of a convergent power series ring over a field $k$ is called an {\it analytic} ring over $k$.
Any complete local ring containing a field is an analytic ring over its coefficient field, and it is known that any analytic local ring is Henselian; see \cite[Chapter VII]{Nagata}.
Yoshino and Kawamoto observed the decomposability of the fundamental module of an analytic normal domain.

\begin{thm}[Yoshino-Kawamoto]
Let $R$ be an analytic normal local domain of dimension two.
Suppose that the residue class field of $R$ is algebraically closed and has characteristic zero.
Then the following conditions are equivalent:
\begin{enumerate}
\item[$(1)$]\ The fundamental module of $R$ is decomposable;
\item[$(2)$]\ $R$ is an invariant subring of a regular local ring by a cyclic group.
(In other words, $R$ is a cyclic quotient singularity.)
\end{enumerate}
\end{thm}

For the details of this theorem, see \cite[Theorem (2.1)]{YK} or \cite[Theorem (11.12)]{Yoshino}.
With the notation of the above theorem, suppose in addition that $R$ is a complete Gorenstein ring such that $\Omega _R ^2 k$ is decomposable.
Then it is seen from Proposition \ref{fund}(4) that $R$ satisfies the condition (1) in the above theorem.
Hence the proof of the above theorem shows that $R$ is of finite Cohen-Macaulay type; see \cite{YK} or \cite{Yoshino}.
It follows from a theorem of Herzog \cite{Herzog} that $R$ is a hypersurface.
Therefore the local ring $R$ is a rational double point of type $(A_n)$ for some $n\geq 1$ by \cite[Proposition (4.1)]{YK}, namely, $R\cong k[[X, Y, Z]]/(XY-Z^{n+1})$.
Thus, the second condition of Theorem \ref{2dimgor} holds.

\begin{ac}
The author would like to give hearty thanks to Shiro Goto.
He motivated the author to consider such a ring as in Proposition \ref{ideal}, which is the dawn of this paper.
The author also wishes to express his deep gratitude to Yuji Yoshino, who gave the author a lot of valuable comments.
In particular, Claim \ref{mn} in the proof of Theorem \ref{2dimgor} is due to him.
The author would also like to thank an anonymous referee very much for useful comments and suggestions.
\end{ac}

%%%%%%%%%%%%%%%%%%% References %%%%%%%%%%%%%%%%%%%%

\end{document}